\documentclass[conference]{IEEEtran}
\IEEEoverridecommandlockouts
\usepackage{cite}
\usepackage{amsmath,amssymb,amsfonts}
\usepackage{algorithmic}
\usepackage{graphicx}
\usepackage{textcomp}
\usepackage{xcolor}

\def\BibTeX{{\rm B\kern-.05em{\sc i\kern-.025em b}\kern-.08em
    T\kern-.1667em\lower.7ex\hbox{E}\kern-.125emX}}
\begin{document}

\title{Static Output Feedback in an Anisotropic-Norm Setup Revisited}

\author{\IEEEauthorblockN{Adrian-Mihail Stoica}
\IEEEauthorblockA{\textit{Faculty of Aerospace Engineering} \\
\textit{University "Politehnica" of Bucharest}\\
Bucharest, Romania \\
adrian.stoica@upb.ro}
\and
\IEEEauthorblockN{Isaac Yaesh}
\IEEEauthorblockA{\textit{Control Dept.} \\
\textit{Elbit Systems}\\
Ramat Hasharon, Israel \\
itzhak.yaesh@elbitsystems.com}

}

\maketitle

\begin{abstract}
The synthesis problem of static output feedback controllers within the anistropic-norm setup is revisited. A tractable synthesis approach involving iterations over a convex optimization problem is suggested, similarly to existing results for the $H_\infty$-norm minimization case. The results are formulated by a couple of Linear Matrix Inequalities coupled via a bilinear equality, revealing, as in the $H_\infty$ case the duality of between the control-type and filtering type LMIs and allowing a tractable iterative method to cope with practical static output feedback synthesis problems. The resulting optimization scheme is then applied to a flight control problem, where the merit of the anisotropic norm setup is shown to provide a useful trade-off between closed loop response and feedback gains.
\end{abstract}

\begin{IEEEkeywords}
Static Output Feedback, $H_\infty$, Anisotropic Norm.  
\end{IEEEkeywords}

\section{Introduction}
The problems of optimal control and filtering received much
attention over the years. Solutions for these
problems were presented by Kwakernaak and Sivan \cite{kwak}. Modelling errors were considered in \cite{robust1}. When the external input signals are of white
noise type, $H_2-$norm minimization is applied, leading to
the Kalman filter \cite{Kalman} and Linear Quadratic Gaussian
(LQG) control. An alternative modelling of the exogenous inputs is
based on deterministic bounded energy signals associated with the $H_\infty$-norm based framework (\cite{zames}) applicable 
both to filtering (\cite{Grimble}, \cite{Simon}) and control (\cite{yaeshShaked}). Since many practical applications require an intermediate solution between $H_2$ and $H_\infty$
Since $H_2$ is not entirely suitable when signals are strongly coloured and $H_\infty$ may result in poor performance when these signals are weakly coloured (e.g. white noise), mixed $H_2/H_\infty$ norm minimization becomes useful (see, e.g. \cite{bernstein}, \cite{rotstein}). A promising alternative to accomplish such
compromise is to use the so-called \emph{a-anisotropic norm} (\cite{kurd1}, \cite{kurd3},\cite{kurd10}) offering
an intermediate topology between the $H_2$ and $H_\infty$ norms.
More precisely, consider the $m$ dimensional coloured signal $w(t), t=0,1,... $ generated by 
the discrete-time stable filter $G$
\begin{eqnarray}
\label{e1}
\begin{array}{rcl}
x_f(t+1)&=&A_fx_f(t)+B_fv(t)\\
w(t)&=&C_fx_f(t)+D_fv(t), t=0,1,...
\end{array}
\end{eqnarray}
where $A_f\in {\mathcal{R}}^{n_f\times n_f}$, $B_f\in {\mathcal{R}}^{n_f\times m}$, $C_f\in {\mathcal{R}}^{m\times n_f}$, $D_f\in {\mathcal{R}}^{m\times m}$  and where $v\in \mathcal{R}$ are independent Gaussian white noises with $E[v(t)]=0$ and $E[v(t)v^T(t)]=I_m$.
Then, the  \emph{a}-anisotropic norm
$|||F|||_a$ of a discrete-time stable system $F$ with the state-space realisation
\begin{eqnarray}
\label{e2}
\begin{array}{rcl}
x(t+1)&=&Ax(t)+Bu(t)\\
y(t)&=&Cx(t)+Du(t), t=0,1,...
\end{array}
\end{eqnarray}
is defined as
\begin{eqnarray}
\label{e4def}
|||F|||_a=\sup_{G\in{\mathcal{G}}_a}\frac{\|FG\|_2}{\|G\|_2},
\end{eqnarray}
 ${\mathcal{G}}_a$ denoting the set of all stochastic systems
of form (\ref{e1}) with the \emph{mean anisotropy} $\bar{A}(G)\leq a$. The mean anisotropy of stationary Gaussian sequences was introduced in \cite{kurd1} and it represents an entropy theoretic measure of the deviation of a probability distribution from Gaussian distributions with zero mean and scalar covariance matrices. In \cite{kurd2}, it is proved based on the Szeg{\"o}-Kolmogorov theorem (\cite{Rozanov}) that the mean anisotropy of a signal generated by an $m$-dimensional Gaussian white noise $v(t)$ with zero mean and identity covariance applied to a stable linear system $G$ with $m$ outputs has the form  
\begin{equation}
\label{andef}
\bar{A}(G) ={ -\frac{1}{2} \ln \,
\det\left(\frac{mE\left[\tilde w(0) \tilde w(0)^T\right]}{Tr\left(E\left[w(0) w(0)^T\right]\right)}\right)},
\end{equation}
where  $E[{\tilde{w}}(0){\tilde{w}}(0)^T]$ is the covariance of the prediction error ${\tilde{w}}(0):=w(0)-E[w(0)|(w(k), k<0]$. In the case when the output $w$ of the filter $G$ is a zero mean Gaussian white noise (i.e. its optimal estimate is just zero), $w(0)$ cannot be estimated from its past values and ${\tilde{w}}(0)=w(0)$ which leads to $\bar{A}(G)=0$.

It is proved (see, for
instance \cite{kurd3}) that the anisotropic norm has the property:
\begin{equation}
\label{e00}
\frac{1}{\sqrt{m}}||F||_2=|||F|||_0 \leq |||F|||_a \leq
||F||_\infty = \lim_{a \rightarrow \infty} |||F|||_a
\end{equation}
In \cite{kurd2}, conditions for the anisotropic norm boundedness are given in terms of a non convex optimization problem while in \cite{kurd10} a convex form of the Bounded Real Lemma (BRL) type result with respect to the anisotropic norm was obtained. One of the leading motivations to use the anisotropic norm is the fact $|||F|||_a
\leq ||F||_\infty$ making it a relaxed version of the
$H_\infty$-norm for many practical cases in which the driving noise
signals can be characterised not just by their finite energy, but
as outputs of a colouring linear systems in a certain class, where
the colouring filters are of a finite anisotropy. In a case study presented in  \cite{Tchaik} it is shown that for a TU-154 type aircraft landing
system,  the $H_\infty$ controller is more efficient than the
corresponding $H_2$ controller for a windshear profile (which is a
coloured rather than a white noise process) but, as could be
expected, is more conservative, in the sense of higher gains and
subsequently larger control actions; moreover,   the anisotropic-norm based controller (based
on an appropriate anisotropic norm bound) is less conservative
than the $H_\infty$ controller and requires significantly smaller
control actions.

The aim of the present paper is to derive a tractable  characterisation
 of the static output feedback synthesis problem. Static Output Feedback (SOF) synthesis is very useful, when the design of fixed structured controllers such as PID (Proportional, Integral, Derivative) is required. Such controllers are very common both in the process control and aerospace control applications and have been adapted to even cope with flexible modes  \cite{icinco2013}. However, it is well known that SOF synthesis problem being non-convex can not be efficently solved and is even an NP-hard problem in some cases \cite{blondel}. Nevertheless, under the $H_2$ and $H_\infty$ setups, the software package MATLAB provides efficient solutions such as $hinfstruct$ procedure \cite{apkarian} and HIFOO Toolbox (\cite{BHLO}). Although the SOF problem within the anisotropic-norm setup has already been considered in \cite{Tchaiks}, we provide a somewhat different solution by a couple of Linear Matrix Inequalities, which although coupled as \cite{Tchaiks} by a bilinear equality, reveal, as in the $H_\infty$ case \cite{iwasaki1}, the duality of between the control-type and filtering type LMIs allowing direct application of the iterative method of \cite{leibfritz} to deal with such a bilinear equality . 
 
 We, therefore, consider the following plant 
 \begin{eqnarray}
 \label{e5bis}
 \begin{array}{rcl}
 x(t+1) &=& Ax(t) + B_1 w(t) + B_2 u(t)
 \end{array}
 \end{eqnarray}
 where we seek for a static control matrix $K$, so that $u(t)=Ky(t)$, where
 $$
 y(t) = C_2 x(t)
 $$
 will minimize
 \begin{eqnarray}
 \begin{array}{rcl}
 \label{e5bisbis}
 z(t)& =& C_1x(t) + D_{12} u(t) + D_{11}w(t)
 \end{array}
 \end{eqnarray}
 in the sense of bounded anisotropic norm, which is yet to be specified.
 
{\textit{Notation.}} Throughout the paper the superscript `$T$'
stands for matrix transposition, ${\mathcal{R}}$ denotes the set
of scalar real numbers whereas ${\mathcal{Z}_+}$ stands for the
non-negative integers. Moreover, ${\mathcal{R}}^n$ denotes the $n$
dimensional Euclidean space, ${\mathcal{R}}^{n\times m}$ is the
set of all $n\times m$ real matrices, and the notation $P\!>\!0$
($P\geq 0$), for $P\in{\mathcal R}^{n\times n}$  means that  $P$
is symmetric and positive  definite (positive semi-definite). The
trace of a matrix $Z$ is denoted by $Tr(Z)$, and $|v|$ denotes
the Euclidian norm of an $n$-dimensional vector $v$. Finally note
that the terms Lyapunov and Riccati equations in this paper, refer
to generalised versions of the standard equations appearing in the
$H_2 $ and $H_\infty$ control literature.

We next state in Section II the problem and provide some preliminaries. An approximate solution is suggested in Section III, whereas an exact but iterative solution is offered in Section IV. Finally,  Section V includes some concluding remarks.

\section{Preliminaries and Motivation}

Consider the following discrete-time system $F$,
described by
\begin{eqnarray}
\label{e66}
\begin{array}{ccl}
x(t+1)&=&Ax(t)+B w(t)\\
y(t)&=& C x(t)  D w(t),\, t=0,1,...
\end{array}
\end{eqnarray}
In the following, some known useful results are briefly reminded.

{\bf Definition 1:}
\label{d2} The $H_2$-type norm of the ES system (\ref{e66}) is
defined as
\begin{eqnarray*}
\|F\|_2=\left[\lim_{\ell\rightarrow\infty}\frac{1}{\ell}\sum_{t=0}^\ell
E\left[y^T(t)y(t)\right]\right]^\frac{1}{2},
\end{eqnarray*}
where $\lbrace y(t)\rbrace_{t\in{\mathcal{Z}}_+}$ is the output of the system
(\ref{e1}) with zero initial conditions generated by the sequence
$\lbrace w(t)\rbrace_{t\in{\mathcal{Z}}_+}$ of independent random vectors with
the property that $E\left[w(t)\right]=0$ and
$E\left[w(t)w^T(t)\right]=I_m$, $\lbrace w(t)\rbrace_{{t\in{\mathcal{Z}}_+}}$.

The next result provides a method to compute the $H_2$ norm of the
system of (\ref{e66}) (see e.g. (\cite{DMS}).

{\bf Lemma 1:}
\label{lem1} The $H_2$ type norm of the ES system (\ref{e6}) is
given by
$\|F\|_2=\left(Tr\left(B^T X B+D^TD\right)\right)^\frac{1}{2}$
where $X\geq 0$ is the solution of the generalised Lyapunov equation
$X=A^T X A + C^T C$.

{\bf Definition 2:} The $H_\infty$ norm of the stable discrete-time system of form (\ref{e66}) is defined as 
\begin{eqnarray*}
\|F\|_\infty=\sup_{\theta\in[0,2\pi)}\lambda_{\max}^{\frac{1}{2}}\left(F^T\left(e^{-j\theta}\right)F\left(e^{j\theta}\right)\right),
\end{eqnarray*}
where $\lambda_{\max}$ denotes the maximal eigenvalue and $F(\cdot)$ is the transfer function of the system.

The $H_\infty$ norm is characterised by the following result, well-known as the Bounded Real Lemma (BRL).

{\bf Lemma 2:}
The stable system (\ref{e66}) has the norm
$\|F\|_\infty<\gamma$ for a certain $\gamma>0$ if and only if the
Riccati equation
\begin{eqnarray*}
\begin{array}{c}
P\!= A^T P A +\left (A^T P B+C^T D \right) \\
\!\!\times\left(\gamma^2I-B^T P B -D^TD\right)^{-1}\\
 \left(A^T P B+C^TD \right)^T+C^TC
\end{array}
\end{eqnarray*}
has a stabilizing solution $P\geq 0$ such that
$\Psi_{1/\gamma^2}:=\gamma^2I-B^T PB-D^TD>0$.

It is recalled (\cite{DMS}) that a symmetric solution $P$ of the
above Riccati equation is called a \emph{stabilising solution} if
the system
\begin{eqnarray*}
x(t+1)=\left( A+BK \right)x(t)
\end{eqnarray*}
is stable, where by definition
\begin{eqnarray*}
K:=\Psi_{1/\gamma^2}^{-1}\left(A^TPB+C^TD\right)^T.
\end{eqnarray*}

To conclude this section, we state the BRL-like result to characterize the anisotropic norm (\cite{kurd2}). Note that for $a$  tending to infinity, the result of Lemma 2 is recovered.

{\bf Theorem 3:}
\label{t1} The system of (\ref{e1}) satisfies $|||F|||_a \leq \gamma$ for a given
$\gamma >0$ if and only if there exists $q\in
\left(0, \min\left(\gamma^{-2},\right.\right.$
$\left.\left.\|F\|_\infty^{-2}\right)\right)$ such that the
Riccati equation
\begin{eqnarray}
\label{e5}
\begin{array}{l}
X=A^T X A + \left(A^T X B +C^TD\right)\\
\times\left(\frac{1}{q}I-B^T X B-D^TD\right)^{-1} \\
\times \left(A^T X B+C^TD\right)^T+C^TC
\end{array}
\end{eqnarray}
has a stabilising solution $X\geq 0$ satisfying the following
conditions
\begin{eqnarray}
\label{e6} \Psi_q:=\frac{1}{q}I-B^TXB-D^TD>0
\end{eqnarray}
and
\begin{eqnarray}
\label{e7} \det\left(\frac{1}{q}-\gamma^2\right)\Psi_q^{-1}\leq
e^{-2a}.
\end{eqnarray}
We conclude this section by providing motivation to use the anisotropic norm. To this end, we denote $\eta =\sqrt{1/q}$ and restate the result of Theorem 3 above, as $||F||_\infty < \eta$
so that $det(\eta^2 - \gamma^2)\Psi_q^{-1}\leq
e^{-2a}$. Namely, $\eta^2-\gamma^2 \leq (det \Psi_q)^{1/m} e^{-2a/m}$. 
Using the general inequality (see \cite{dbern}) $(det \Psi_q)^{1/m} \leq \frac{Tr{\Psi_q}}{m}$ valid for any $\Psi_q\geq 0$, and noting that $\Psi_q = \eta^2 I-B^TXB-D^TD$, the following motivating result of \cite{{StoicaYaesh2020}} was obtained:

{\bf Lemma 4:}
Consider the system $F$ of (\ref{e1}). Let $\eta$ and $\sigma$ respectively satisfy 
$$
||F||_\infty < \eta \mbox{ and } ||F||_2 < \sigma 
$$
The $a$-anisotropic norm of the system of (\ref{e1}) is then upper bounded by the following linear interpolation between its $H_\infty$ and $H_2$ norms. Namely,
$$
\gamma^2 \geq \eta^2(1-e^{-2a/m}) + \frac{\sigma^2}{m} e^{-2a/m}
$$
We, therefore, see that in view of (\ref{e00})\, one may interpret the $a$-anisotropic norm, the following approximate relation
$$
|||F|||_a^2\approx |||F|||_\infty^2 (1-e^{-2a/m}) + |||F|||_0^2  e^{-2a/m}
$$
providing a useful insight to the $a$-anisotropic norm, which approximation can be interpreted also as just mixed $H_\infty /H_2$ optimization, however, in the exact proportions dictated by the Lemma, in terms of $e^{-2a/m}$. 

\section{Static Output Feedback}

 We next provide a solution to the problem of synthesis of SOF control synthesis under the anisotropic norm, using iterative solution for LMIs. To this end, we define the cost function to be
\begin{equation}
    \label{costj}
    J(K) = |||\mathcal{F}_{c\ell}(K)|||_a
\end{equation}
where $\mathcal{F}_{c\ell}(K)$ denotes the closed loop system obtained from (\ref{e5bis}) and (\ref{e5bisbis}) with the static output feedback $u_2(t)=Ky(t)$ having the realization 
\begin{eqnarray*}
\begin{array}{rcl}
x(t+1)&=&(A+B_2KC_2)x(t)+B_1w(t)\\
z(t)&=&(C_1+D_{12}KC_2)x(t)+D_{11}w(t)
\end{array}
\end{eqnarray*}
   Using Theorem 3, if follows that the above closed loop system $\mathcal{F}_{c\ell}(K)$  is stable and it has the $a$-anisotropic norm less than a given $\gamma>0$ if and only if there exist a $q\in\left(0, \min(\gamma^{-2}, \|F_{c\ell}\|_\infty^{-2})\right)$ and a symmetric matrix $X>0$ such that
\begin{eqnarray}
\label{e4}
\left[ \begin{array}{cc}{\mathcal E}_1(X,K)\ & {\mathcal E}_2(X,K) \\
(1,2)^T & -\frac{1}{q}I+B_1^T X B_1+D_{11}^TD_{11}
\end{array}\right]<0
\end{eqnarray}
where
\begin{eqnarray*}
{\mathcal E}_1(X,K):=-X+\left(A+B_2KC_2\right)^TX\left(A+B_2KC_2\right)\\+\left(C_1+D_{12}KC_2\right)^T\left(C_1+D_{12}KC_2\right)
\end{eqnarray*}
\begin{eqnarray*}
{\mathcal E}_2(X,K):=\left(A+B_2KC_2\right)^TXB_1\\+\left(C_1+D_{12}KC_2\right)^TD_{11}
\end{eqnarray*}
and 
\begin{eqnarray}
\label{e5add}
\frac{1}{q}- \gamma^2<e^{-\frac{2a}{m}}\left(\det \left(\frac{1}{q}I-B_{1}^T  B_{1}-D_{11}^TD_{11}\right)\right)^{\frac{1}{m}}
\end{eqnarray}
Based on Schur complements arguments, one can see that the inequality (\ref{e4}) is equivalent with 
\begin{eqnarray}
\label{e6add}
\left[ \begin{array}{cccc}
-X & 0 & \left(3,1\right)^T & \left(4,1\right)^T\\
0 & -\frac{1}{q}I & B_1^T & D_{11}^T\\A+B_2KC_2& B_1 & -X^{-1}& 0\\
C_1+D_{12}KC_2 & D_{11} &0 &-I\end{array}\right]<0
\end{eqnarray}


Multiplying the above inequality to the left and to the right by $diag(I,I,X,I)$ one obtains that is is equivalent with
\begin{eqnarray}
\label{e44}
\left[ \begin{array}{cccc}
-X & 0 & \left(3,1\right)^T & \left(4,1\right)^T\\
0 & -\frac{1}{q}I & B_1^T X& D_{11}^T\\X(A+B_2KC_2)& XB_1 & -X& 0\\
C_1+D_{12}KC_2 & D_{11} &0 &-I\end{array}\right]<0
\end{eqnarray}

which may be re-written as
\begin{equation}
\label{e45}
\mathcal{Z}+{\mathcal{P}}^T K\mathcal{Q}+{\mathcal{Q}}^T K^T\mathcal{P}<0,
\end{equation}
where we denoted
\begin{eqnarray}
\label{e46}
\begin{array}{l}
\mathcal{Z}:=\left[\begin{array}{cccc} -X&  0 &A^TX & C_1^T\\
0&-\frac{1}{q}I &B_1^TX &D_{11}^T\\
XA& XB_1 &-X &0\\
C_1 &D_{11} & 0 &-I\end{array}\right], \\
 {\mathcal{P}}^T:=\left[\begin{array}{c}0\\0\\XB_2\\D_{12}\end{array}\right],\, 
  {\mathcal{Q}}^T:=\left[ \begin{array}{c}C_2^T \\0\\ 0 \\0 \end{array}\right].
  \end{array}
\end{eqnarray}
According with the so-called Projection lemma (see e.g. \cite{SYG}), the inequality (\ref{e45}) is feasible with respect to $K$ if and only if the following conditions are accomplished
\begin{eqnarray}
\label{e47}
W_{\mathcal{P}}^TZ W_{\mathcal{P}}<0
\end{eqnarray}
and
\begin{eqnarray}
\label{e48}
W_{\mathcal{Q}}^TZ W_{\mathcal{Q}}<0,
\end{eqnarray}
where $W_{\mathcal{P}}$ and $W_{\mathcal{Q}}$ are any bases of the null spaces of ${\mathcal{P}}$ and ${\mathcal{Q}}$, respectively. Since a base of the null space of ${\mathcal{P}}$ is :
\begin{eqnarray}
\label{e481}
W_{\mathcal{P}} = \left[ \begin{array}{ccc} I & 0 & 0 \\0 & I & 0\\0 & 0 & X^{-1}W_1 \\0 & 0 & W_2 \end{array} \right] 
\end{eqnarray}
where $W:=\left[ \begin{array}{c} W_1 \\ W_2 \end{array} \right]$ is the orthogonal complement of $\left[ \begin{array}{cc} B_2^T & D_{12}^T \end{array} \right]$.
Similarly, a base of the null space of ${\mathcal{Q}}$ is:
\begin{eqnarray}
\label{e484}
W_{\mathcal{Q}} = \left[ \begin{array}{ccc} W_3 & 0 & 0 \\W_4 & 0 & 0 \\0 & I & 0 \\0 & 0 & I \end{array} \right]
\end{eqnarray}
where $V:=\left[ \begin{array}{c} W_3 \\ W_4 \end{array} \right]$ is the orthogonal complement of $\left[ \begin{array}{cc} C_2 & 0 \end{array} \right]$.
In order to simplify the inequality of (\ref{e47}) we next express 
$$
W_{\mathcal{P}} = \left[ \begin{array}{cccc} I & 0 & 0 & 0 \\0 & I & 0 & 0 \\0 & 0 & W_1^T & W_2^T \end{array} \right] \left[ \begin{array}{cccc} I & 0 & 0 & 0 \\0 & I & 0 & 0 \\0 & 0 & X^{-1} & 0 \\0 & 0 & 0 & I \end{array}\right] 
$$
Namely 
$$
W_{\mathcal{P}} =\left[ \begin{array}{cc} I & 0 \\0 & W^T  \end{array} \right] \left[ \begin{array}{cccc} I & 0 & 0 & 0 \\0 & I & 0 & 0 \\0 & 0 & X^{-1} & 0 \\0 & 0 & 0 & I \end{array}\right]
$$
with the above definition of $W$. 
Therefore, (\ref{e47}) is simply expressed as
\begin{eqnarray}
\label{e21bis}
\left[ \begin{array}{cc} I & 0 \\ 0 & W^T \end{array} \right] \left[ \begin{array}{cc} \mathcal{M}_{11} & \mathcal{M}_{12} \\ \mathcal{M}_{12}^T & \mathcal{M}_{22} \end{array} \right] \left[ \begin{array}{cc} I & 0 \\ 0 & W \end{array} \right] < 0
\end{eqnarray}
where $\mathcal{M}$ is given by  
$$
\mathcal{M}:=\left[ \begin{array}{cc} \mathcal{M}_{11} & \mathcal{M}_{12} \\ \mathcal{M}_{12}^T & \mathcal{M}_{22} \end{array} \right]
$$
with 
\begin{eqnarray}
\label{e488}
\begin{array}{l}
\mathcal{M}_{11} = \left[ \begin{array}{cc} -X & 0 \\0 &  -\frac{1}{q}I \end{array} \right] ,
\\ \mathcal{M}_{12} = \left[ \begin{array}{cc} A & B_1 \\C_1 &  D_{11} \end{array} \right]^T , \mathcal{M}_{22} = \left[ \begin{array}{cc} -X^{-1} & 0 \\0 &  -I \end{array} \right]
\end{array}
\end{eqnarray}
Then, from (\ref{e21bis}), using the Schur complement of ${\mathcal{M}}_{11}$ it follows that
or by using Schur complements
\begin{eqnarray}
\label{e22bis}
W^T (\mathcal{M}_{22} - \mathcal{M}_{12}^T \mathcal{M}_{11}^{-1} \mathcal{M}_{12} )W <0
\end{eqnarray}
Substituting the definition for $\mathcal{M}_{ij}, i,j=1,2$ and recalling the definition $\eta^2 = \frac{1}{q}$ we obtain the following convenient form of (\ref{e22bis}) which we denote by ${\mathcal L}_Y<0$ the inequality
$$
W^T\left[ \begin{array}{cc} -Y + A Y A^T + B_1 B_1^T & A Y C_1^T + B_1 D_{11}^T \\C_1 Y A^T + D_{11} B_1^T & -\Phi_Y \end{array} \right]W<0
$$
where
$$
\Phi_Y:=\eta^2 I-C_1 Y C_1^T-D_{11} D_{11}^T 
$$
and where we have defined
\begin{eqnarray}
\label{exy1}
\eta^{-2}Y=X^{-1}
\end{eqnarray}
We next repeat the same lines to simplify (\ref{e48}) as well. To this end, we partition
$$
W_Q = \left[ \begin{array}{cc} V & 0 \\ 0 & I \end{array}\right]  
$$
and readily obtain using Schur complements, that (\ref{e48}) is equivalent to 
$$
V^T(\mathcal{N}_{11}-\mathcal{N}_{12}\mathcal{N}_{22}^{-1}\mathcal{N}_{12}^T)V<0
$$
where
\begin{eqnarray}
\label{e499}
\begin{array}{l}
\mathcal{N}_{11} = \left[ \begin{array}{cc} -X & 0 \\0 &  -\frac{1}{q}I \end{array} \right] ,
 \mathcal{N}_{12} = \left[ \begin{array}{cc} XA & XB_1 \\C_1 &  D_{11} \end{array} \right]^T , \\
 \mathcal{N}_{22} = \left[ \begin{array}{cc} -X & 0 \\0 &  -I \end{array} \right]
\end{array}
\end{eqnarray}
We, therefore, obtain the following form of (\ref{e48}) 
$$
V^T\left[ \begin{array}{cc} -X + A^T X A + C_1^T C_1 & A^T X B_1 + C_1^T D_{11} \\B_1^T X A^ + D_{11}^T C_1 & -\Phi_X \end{array} \right]V<0
$$
where
$$
\Phi_X:=\eta^2 I-B_1^T X B_1 -D_{11}^T D_{11}
$$
We denote the latter inequality by :
$$
{\mathcal L}_X < 0 .
$$
We summarize the above results in the following result.

{\bf Theorem 5:} The closed loop system ${\mathcal{F}}_{c\ell}(K)$ is stable and it has the $a$-anisotropic norm less than a given $\gamma>0$ if there exist  symmetric matrices $X>0$ and $Y>0$ and a scalar $\eta^2$ which satisfy the following dual LMIs 
$$
{\mathcal L}_X < 0 , {\mathcal L}_Y < 0 
$$
so both the following convex condition 
\begin{eqnarray}
\label{e25bis}
\eta^2 - det(\Phi_X)^{1/m}e^{-2a/m}< \gamma^2
\end{eqnarray}
and the additional bilinear condition :
\begin{equation}
\label{exy}
XY = \eta^2 I 
\end{equation}
are satisfied.

If the conditions of Theorem 5 are satisfied then the static output gain  may be obtained solving (\ref{e45}) with respect to $K$.

However, the above requires a solution of a set of Bilinear Matrix Inequalities (BMI) due to the $XY=\eta^2 I$ equality. One way to tackle the BMI is to adopt a by first relaxing $XY=\eta^2 I$ by 
\begin{eqnarray}
\label{e8}
\left[ \begin{array}{cc}
X &\eta I\\
\eta I & Y\end{array}\right]>0
\end{eqnarray}
Then if one minimizes $Tr\{XY\}$, the bilinear constraint of (\ref{exy}) is satisfied. To this end, a sequential linearization algorithm (see e.g. \cite{peaucelle1}) can be used. In the initialization step,  the convex problem comprised of the inequalities (\ref{e7}), (\ref{e8}) and (\ref{e5}) is solved for a given $\gamma >0$, and $k=0$, $X_k = 0$ and $Y_k = 0$ are set. Next step where $k$ is set to $k+1$ and $X,Y$ are found so as to minimize
$$
f_k:=Tr\{X_k Y + X Y_k\}
$$
subject to $\mathcal{L}_X<0$, $\mathcal{L}_Y<0$, (\ref{e25bis}) and (\ref{e8}). Then $X_k = X$ and $Y_k = Y$ are set. This step is repeated until $f_k$ is small enough. A related algorithm requiring also line search but with improved convergence properties has been suggested in \cite{leibfritz} and will be applied in the calculations below of the numerical example of the next section. To this end we note that one could choose also $X = \tilde Y^{-1}$ rather than $\eta^{-2}Y=X^{-1} $. In such a case, $\mathcal{L}_Y<<0$ is replaced by
$$
W^T\left[\begin{array}{cc} -\tilde Y + A \tilde Y A^T + q B_1 B_1^T & A \tilde Y C_1^T + q B_1 D_{11}^T \\C_1 Y A^T + q D_{11} B_1^T & -\Phi_{\tilde Y} \end{array} \right]W<0
$$
where $\Phi_{\tilde Y}$ is defined to be
$$
\Phi_{\tilde Y}:=I-C_1 \tilde Y C_1^T-q  D_{11} D_{11}^T 
$$
and (\ref{e8}) should be replaced by
\begin{eqnarray}
\label{e8n}
\left[ \begin{array}{cc}
X & I\\
I & \tilde Y\end{array}\right]>0
\end{eqnarray}
Although this different choice of $Y$ reveals the duality between the control and filtering type inequalities in a less obvious manner, it is more convenient to deal with. Note that to apply \cite{leibfritz} one needs to define a new variable $q=\eta^{-2}$ where in addition to $X \tilde Y=I$ also the scalar valued bilinear equality constraint $\eta^2 q = 1$ has to be satisfied. To this end, one needs also to consider the relaxed version  
\begin{eqnarray}
\label{e9n}
\left[ \begin{array}{cc}
\eta^2  & 1\\
1 & q \end{array}\right]>0
\end{eqnarray}
so that the minimization steps involve now searching for $X,Y, \eta^2, q$ are found so as to minimize
$$
f_k:=Tr\{X_k Y + X Y_k\} + \eta_k^2 q+ \eta^2 q_k\}
$$
subject to $\mathcal{L}_X<0$, $\mathcal{L}_{\tilde Y}<0$, (\ref{e8n}) and (\ref{e9n}).

\section{Application to Flight Control}

We next consider the numerical example of \cite{ackermann}with the synthesis of pitch control loop for the F4E aircraft. Consider 
\[
\begin{array}{l}
 {\rm{                          }}\frac{{\rm{d}}}{{{\rm{dt}}}}\left[ {\begin{array}{*{20}c}
   {{\rm{N}}_{\rm{z}} }  \\
   {\rm{q}}  \\
   {{\rm{\delta }}_{\rm{e}} }  \\
\end{array}} \right] = \left[ {\begin{array}{*{20}c}
   {a_{{\rm{11}}} } & {a_{{\rm{12}}} } & {{\rm{a}}_{{\rm{13}}} }  \\
   {a_{{\rm{21}}} } & {a_{{\rm{22}}} } & {{\rm{a}}_{{\rm{23}}} }  \\
   {\rm{0}} & {\rm{0}} & { - {\rm{30}}}  \\
\end{array}} \right]\left[ {\begin{array}{*{20}c}
   {{\rm{N}}_{\rm{z}} }  \\
   {\rm{q}}  \\
   {{\rm{\delta }}_{\rm{e}} }  \\
\end{array}} \right] +\\
\left[{\begin{array}{*{20}c}
   {{\rm{b}}_{\rm{1}} }  \\
   {\rm{0}}  \\
   {{\rm{30}}}  \\
\end{array}} \right]{\rm{u}} + \left[ {\begin{array}{*{20}c}
   {\rm{1}} & {\rm{0}} & {\rm{0}}  \\
   {\rm{0}} & {\rm{1}} & {\rm{0}}  \\
   {\rm{0}} & {\rm{0}} & {\rm{1}}  \\
\end{array}} \right]{\rm{\omega }}\vspace{.2cm} \\
 {\rm{                               z}} = \left[ {\begin{array}{*{20}c}
   {\rm{1}} & {\rm{0}} & {\rm{0}}  \\
   {\rm{0}} & {\rm{1}} & {\rm{0}}  \\
   {\rm{0}} & {\rm{0}} & {\rm{0}}  \\
\end{array}} \right]{\rm{x}} + \left[ {\begin{array}{*{20}c}
   {\rm{0}}  \\
   {\rm{0}}  \\
   {\rm{0.001}}  \\
\end{array}} \right]{\rm{u}}\vspace{.2cm} \\
 {\rm{                               y}} = \left[ {\begin{array}{*{20}c}
   {\begin{array}{*{20}c}
   {\rm{1}}  \\
   {\rm{0}}  \\
\end{array}} & {\begin{array}{*{20}c}
   {\rm{0}}  \\
   {\rm{1}}  \\
\end{array}} & {\begin{array}{*{20}c}
   {\rm{0}}  \\
   {\rm{0}}  \\
\end{array}}  \\
\end{array}} \right]{\rm{x}} \\
  \end{array}\]

The state-vector consists, of the load-factor $N_z$ ,the pitch-rate
$q$ and elevon angle $\delta_e$. 
The latter relates to the elevon command $u$ via a first-order servo model of a bandwidth of $30 rad/sec$. We note that we have deliberately chosen a small weight on $u$ in $z$ to obtain small damping ratio of the closed-loop poles, when no pole-placement restrictions are imposed. The parameters $a_{i,j}, i=1,2; j=1,2,3, b_1$ are given in
\cite{ackermann} at the four operating points listed in the following table:
\newline
{\small
\begin{tabular}{||c|c|c|c|c||}
\hline\hline
 Operating point & 1   &2  & 3 &4\\
 \hline
 Mach number&.5&.9&.85&1.5\\
 \hline
 Altitude (ft)&5000&35000&5000&35000\\
\hline \hline
$a_{11}$ & -.9896 &-.6607 &-1.702 &-.5162\\
\hline
 $a_{12}$ &
17.41&18.11&50.72&29.96\\
\hline
 $ a_{13}$ & 96.15 &
84.34 & 263.5 & 178.9\\
\hline
 $a_{21}$ & .2648 &.08201 &.2201 &-.6896\\
\hline
 $a_{22}$ & -.8512 &-.6587 &-1.418 &-1.225\\
 \hline
 $a_{23}$ &-11.39
&-10.81 & -31.99 &-30.38\\
\hline
 $b_1$ & -97.78 &-272.2 & -85.09 &-175.6\\
\hline\hline
\end{tabular}}

\vspace{0.3cm}

{\bf Table 7:} The parameters of the four operating points.

\vspace{0.3cm}

The static output controller $u=Ky$ with $K = \left[ \begin{array}{cc} K_1 & K_2 \end{array} \right]$ will be designed for each of the four operating points, applying Theorem 5, applying the iterative procedure by \cite{leibfritz} to deal with the bilinear equality. We will first take mean anisotropy level of $a$ that tends to $\infty$, to obtain the $H_\infty$ controller. Next we take  $a=0.5$. The upper bound on $\gamma$ for $H_\infty$ control was taken as envelope point dependent $\gamma_{points \infty} = \left[ \begin{array}{cccc} 0.3 & 0.6  &  1 & 0.25 \end{array} \right]$. The corresponding upper bounds for the a-anisotropic norm are generally lower, with a point-wise factor of $\left[ \begin{array}{cccc} 0.84 & 0.8 & 0.82 & 1 \end{array} \right]$. The results of those designs are depicted in Fig. 1 and Fig.2 where  closed-loop singular values as well as the design bound and the actual norm ($H_\infty$ and a-anisotropic) are shown. As could be expected, the maximum singular values are somewhat lower for the $H_\infty$ design, at the cost, however, of larger gains with respect to the a-anistropic design. The gains comparison of the Table 8 demonstrates the advantage offered by the a-anisotropic design, at the cost in the closed loop maximum singular values indicated above. However, if the same bounds as those taken for the $H_\infty$ design, are taken for the a-anistropic design, namely, a sub-optimal a-anisotropic case is considered, then closed loop singular values of Fig. 3. are obtained. The gains of those three designs are given in Tables 8, 9 and 10 respectively. To allow an easy comparison, Fig, 4 depicts the singular values of all three designs together. A close scrutiny of the gains and the singular values reveals that the sub-optimal a-anisotropic design allow a fine trade-off between the maximum singular value and the gain values. Namely, considerably lower gain vector norms are obtained at the cost of a rather small sacrifice in the maximum singular value.    

\section{Conclusions}
A synthesis scheme for static output feedback controllers has been derived, under the setup of $a$-anisotropic norm which is based on an intermediate topology between $H_2$ and $H_\infty$. Given a required norm-bound, set of Linear Matrix Inequalities, along with a geometric-mean convex inequality, and an additional bilinear equality, characterize sub-optimal controllers as in\cite{Tchaiks}, however revealing the duality in the style of \cite{iwasaki1} between the control and filtering type Linear Matrix Inequalities. The latter dual form is convenient to use with the iterative algorithm of \cite{leibfritz} to alleviate the bilinear equality.  The resulting control design procedure is useful e.g. in the aerospace industry for flight control loops, where controllers with classical "cook-book" structures in the style of  Proportional-Integral-Derivative controllers. An example from the field of flight control is given comparing $H_\infty$ design with a couple of a-anisotropic designs applying the above mentioned iterative algorithm of \cite{leibfritz}. The results suggest that the latter offers a considerable reduction in the gains ate the cost of small increase in the closed-loop maximum singular values. Those results encourage further experimenting with the a-anisotropic setup for static-output feedback control and motivates further research also along the lines of \cite{Tchaik} and \cite{Tchaiks} and the references therein.

\newpage
\begin{figure}[h]
\centering
\includegraphics[width=8cm,height=8cm]{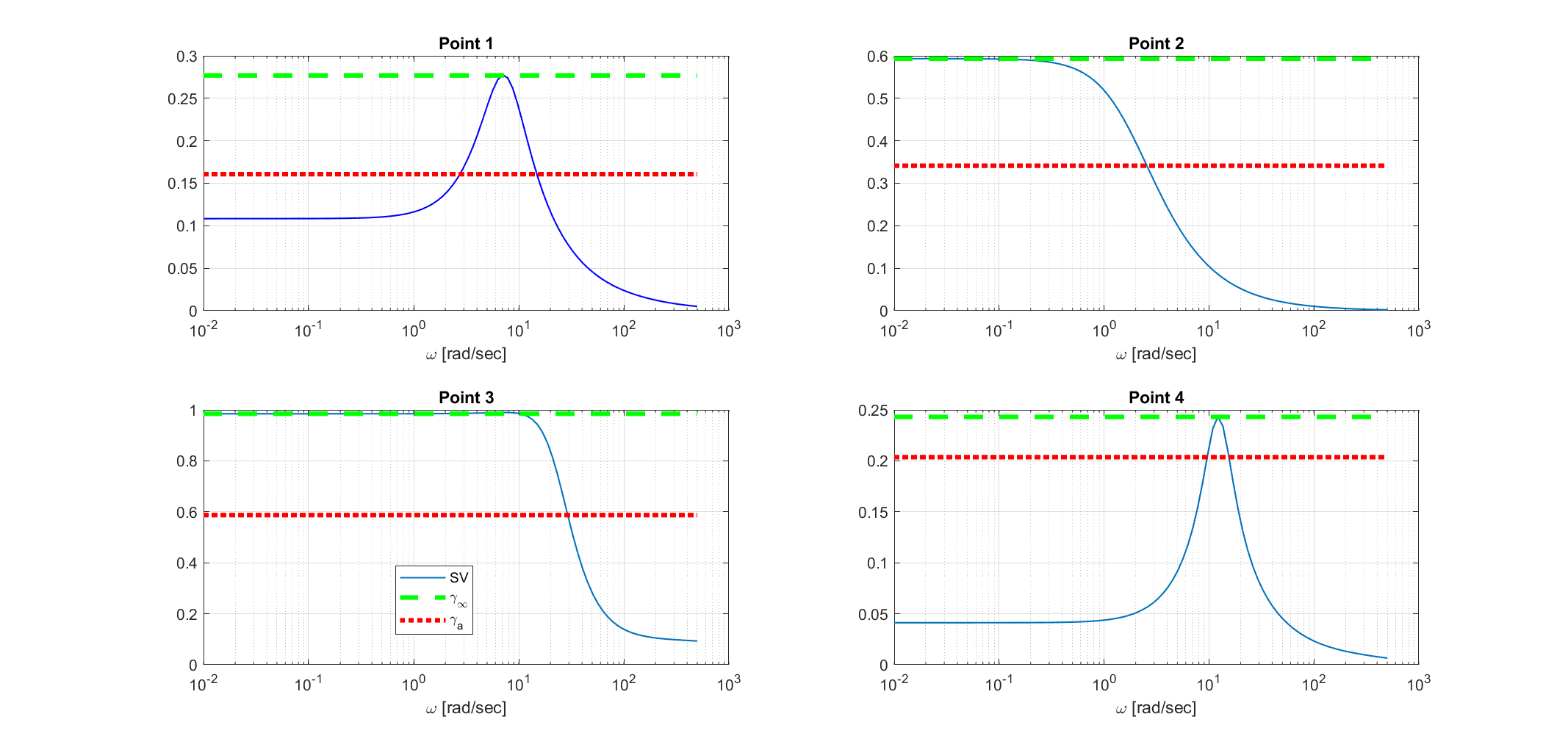}
\caption{$T_{zw} - H_\infty Design$ }
 \end{figure}

\begin{figure}[h]
\centering
\includegraphics[width=8cm,height=8cm]{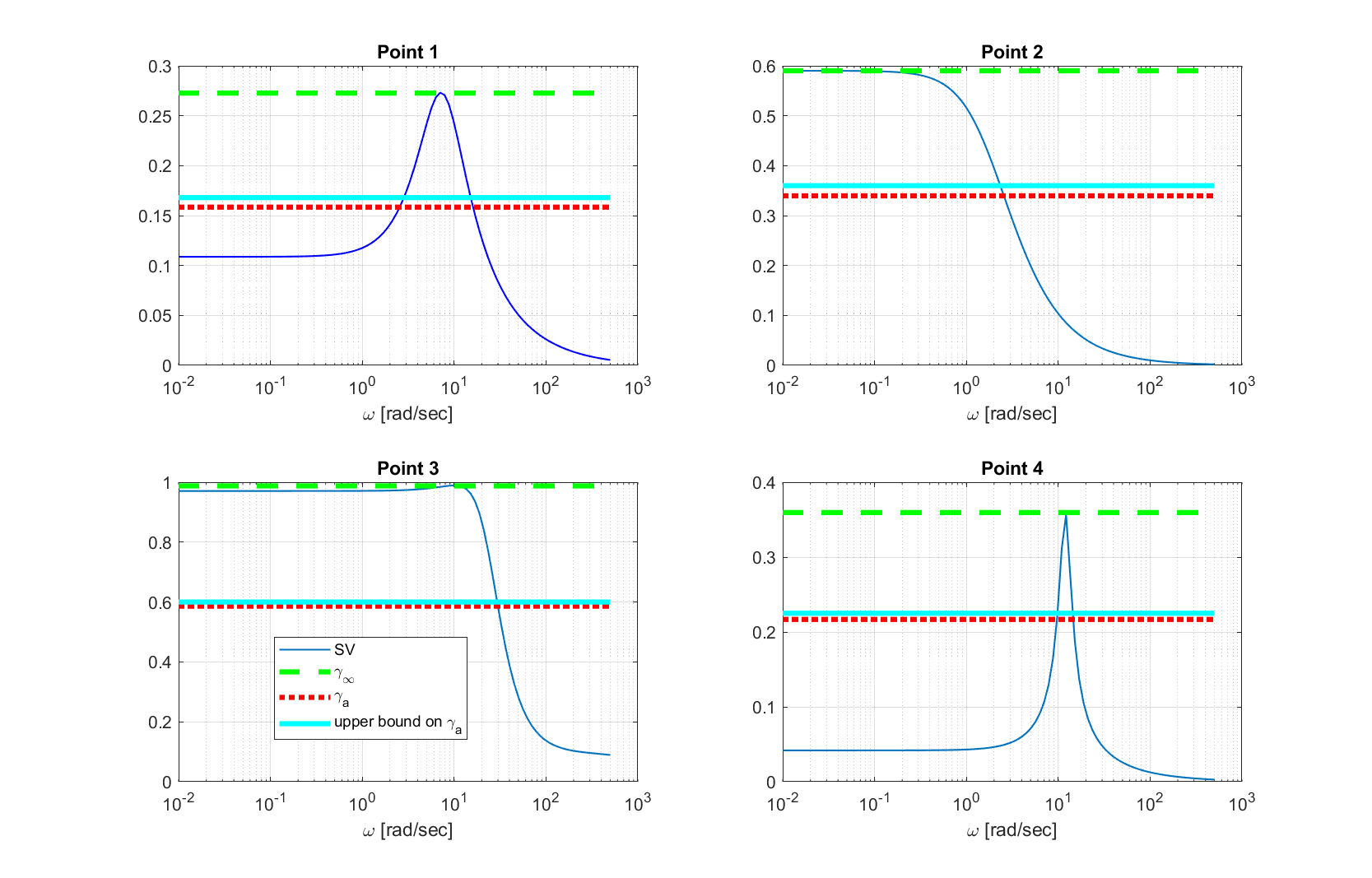}
\caption{$T_{zw}$ Near Optimal- a-Anisotropic Design}
 \end{figure}

\begin{figure}[h]
\centering
\includegraphics[width=8cm,height=8cm]{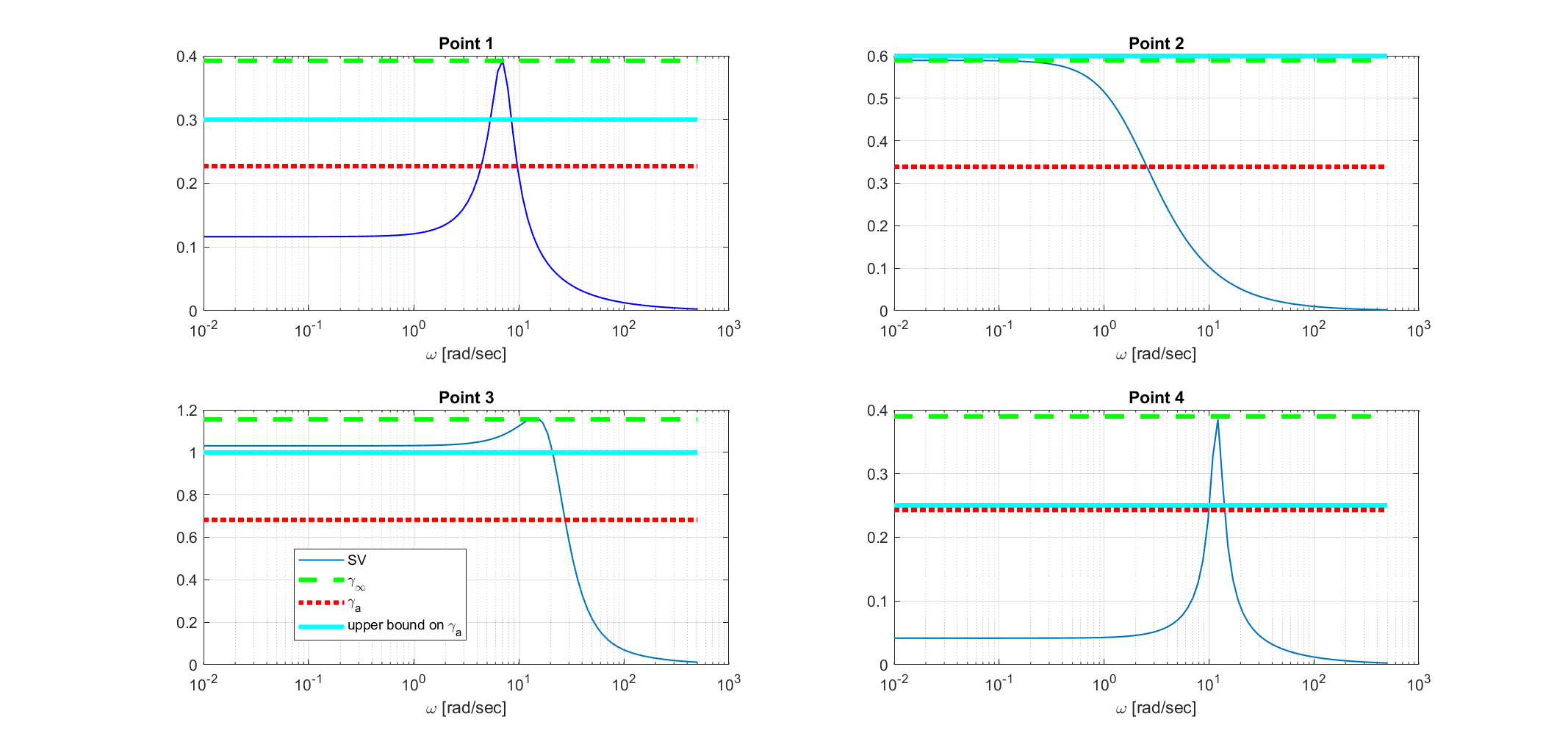}
\caption{$T_{zw}$ Sub-Optimal- a-Anisotropic Design}
 \end{figure}

\begin{figure}[h]
\centering
\includegraphics[width=8cm,height=8cm]{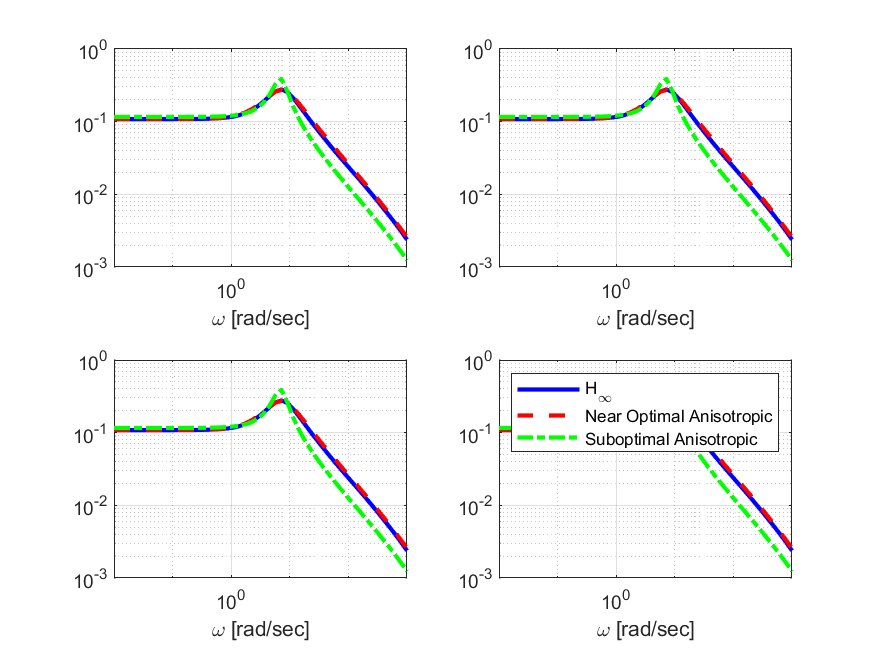}
\caption{$T_{zw}$ Sub-Optimal- a-Anisotropic Design}
 \end{figure}

\vspace{0.3cm}

\end{document}